\documentclass[a4paper,oneside,10pt]{article}%
\usepackage{amsmath}
\usepackage{amsfonts}
\usepackage{amssymb}
\usepackage{graphicx}
\usepackage[square,numbers,sort&compress]{natbib}%
\providecommand{\U}[1]{\protect\rule{.1in}{.1in}}

\newtheorem{theorem}{Theorem}[section]

\newtheorem{assumption}[theorem]{Assumption}

\newtheorem{proposition}[theorem]{Proposition}
\newtheorem{remark}[theorem]{Remark}

\numberwithin{equation}{section}

\begin{document}

\title{\textbf{{Stochastic Linear Quadratic Optimal Control with General Control
Domain}}}
\author{Shaolin Ji \thanks{Zhongtai Securities Institute for Financial Studies, Shandong University, Jinan
250100, China and Institute of Mathematics, Shandong University, Jinan 250100,
China, Email: jsl@sdu.edu.cn. This work was supported by National Natural
Science Foundation of China (No. 11571203); Supported by the Programme of
Introducing Talents of Discipline to Universities of China (No.B12023). }
\and Xiaole Xue\thanks{Corresponding author. Zhongtai Securities Institute for Financial Studies,
Shandong University, Jinan 250100,
China. Email: xiaolexue1989@gmail.com, xuexiaole.good@163.com.  \newline
   \textit{3 Jun 2017 - Submitted to Math. Control Relat. Fields}}}
\maketitle

\textbf{Abstract}. This paper considers the stochastic linear quadratic
optimal control problem in which the control domain is nonconvex. By the functional analysis and convex
perturbation methods, we establish a novel maximum principle. The application of
the proposed maximum principle is illustrated through a work-out example.



\addcontentsline{toc}{section}{\hspace*{1.8em}Abstract}

\section{Introduction}

The stochastic linear
quadratic optimal control problems play an important role in optimal control
problems. On one hand, many nonlinear control problems can be approximated by
the linear quadratic control problems; on the other hand, solutions to the
linear quadratic control problems show elegant properties because of their
brief and beautiful structures.

A classical form of the stochastic linear quadratic optimal control problems
is to minimize the quadratic cost functional
\begin{equation}%
\begin{array}
[c]{rcl}%
J(u(\cdot)) & = & \mathbb{E}\{\frac{1}{2}\int_{0}^{T}%
[<Q(t)X(t),X(t)>+2<S(t)X(t),u(t)>+<R(t)u(t),u(t)>]dt\\
&  & +\frac{1}{2}<GX(T),X(T)>\}
\end{array}
\label{classic obj}%
\end{equation}
with the control $u(\cdot)$ being a square integrable adapted process and the
state being the solution to the linear stochastic differential equation
\begin{equation}
\left\{
\begin{array}
[c]{rcl}%
dX(t) & = & [A(t)X(t)+B(t)u(t)+b(t)]dt+[C(t)X(t)+D(t)u(t)+\sigma
(t)]dW_{t},\ \ t\in\lbrack0,T],\\
X(0) & = & x,
\end{array}
\right.  \label{classic state}%
\end{equation}
where $A,B,C,D,b,\sigma$ are deterministic matrix-valued functions of suitable
sizes; $T$ is a fixed terminal time; $W$ is a Wiener process. Under suitable
conditions, the cost functional $(\ref{classic obj})$ and the control system
$(\ref{classic state})$ are well defined. This kind of stochastic linear
quadratic problem was first studied by Wonhanm \cite{wonham}. Then Bismut
\cite{bismut} proved an existence result for optimal control.

As one of the foundations of modern finance theory,
the mean-variance portfolio issue can be translated to the
stochastic linear quadratic problems. As indicated in
\cite{lizhoulim}, the mean-variance problem is essentially
an indefinite stochastic linear quadratic problem because its running
cost is identically zero.
The mean-variance model was formulated as a quadratic programming problem by
Markowitz in the case of a single-period
investment (see \cite{markowitz1,markowitz2}), and the multi-period
counterpart was soon considered after the pioneering
work (see \cite{elton,grauer,samuelson,mossin}).
By resorting to the method of dynamic programming,
the mean-variance model in a continuous time setting was developed a bit later (see
\cite{follmer,duffie-j,duffie-r}). By employing the embedding technique, Zhou and Li \cite{zhouli}
turned the continuous-time mean-variance problem into a
stochastic linear quadratic problem, which could be
solved by using the results in  Chen, Li and Zhou \cite{chenlizhou}.

Recently, Li, Zhou and Lim \cite{lizhoulim} discussed a continuous-time
mean-variance problem with no-shorting constraints. By using a
martingale approach, Bielecki et al. \cite{bielecki} investigated the problem of
continuous-time mean-variance with bankruptcy prohibition.
Heunis \cite{heunis} carefully considered to minimize the expected value of a general quadratic
loss function of the wealth in a more general setting where there is a
specified convex constraint on the portfolio over the trading interval,
together with a specified almost sure lower-bound on the wealth at close of
trade. Li and Xu \cite{lixu} studied the continuous-time mean-variance problem
with the mixed restriction of both bankruptcy prohibition and convex cone
portfolio constraints. The control domains in the aforementioned literature are assumed to be convex.
Inspired by the above discussion, we intend to further study the stochastic linear quadratic problem
with nonconvex control domain, which may be helpful to the mean-variance
problems with nonconvex control domain.

In this paper, we develop a novel stochastic maximum principle to tackle the
stochastic linear quadratic problem with nonconvex control domain. To begin
with, based on functional analysis approach, we transform the original linear
quadratic problem into a quadratic optimization problem in a Hilbert space,
which consists of all square integrable control processes. Next, by
introducing a parameter in the quadratic optimization problem, we turn
equivalently the original problem into a concave control problem with convex
control domain, which can be solved by the classical stochastic maximum
principle. Then we derive a new stochastic maximum principle. It is well known that the stochastic maximum principle is of great importance in solving stochastic optimal
control problems (see \cite{ElKaroui}, \cite{jizhou}, \cite{peng1990,peng1993}%
, \cite{wu2013,yong2010}). A local form of the stochastic maximum principle for the
classical stochastic recursive optimal control problem was established later in Peng \cite{peng1993}. By utilizing the spike variational method
and the second-order adjoint equations, Peng \cite{peng1990} obtained a general stochastic
maximum principle when the admissible control domain needs not to be convex and the
diffusion coefficients contain the control variable.
The achievement of the stochastic maximum principle in \cite{peng1990}
relies heavily on the second-order adjoint equations. Compared with the classical methods dealing with the nonconvex stochastic linear quadratic
problem, the approach provided in this paper has two main advantages as follows.
Firstly, we need not to introduce the second-order adjoint equation, and the
presented stochastic maximum principle has a concise form. Secondly, we don't
have to impose any positive definite on the coefficients.

The paper is organized as follows. In section 2, we present some notations and
formulate the stochastic linear quadratic problem. We formulate our stochastic
linear quadratic optimal control problem by functional analysis approach in
section 3. We derive the maximum principle in section 4 and give a work-out
example in section 5. In section 6 we conclude the paper.

\section{Problem Formulation}

Throughout this paper, we denote by $\mathbb{R}^{n}$ the $n-$dimensional
vector space and $\mathbb{R}^{k\times n}$ the set of $k\times n$ matrices.
Particularly, we denote by $\mathbb{S}^{n}$ the set of symmetric $n\times n$
matrices. For any given Euclidean space $H$, we denote by $\langle\cdot
,\cdot\rangle$ (resp. $\mid\cdot\mid$) the scalar product (resp. norm) of $H$.
When $M=(m_{ij}),N=(n_{ij})\in\mathbb{R}^{k\times n}$, we define
$<M,N>=tr\{MN^{\intercal}\}$ and $\mid M\mid=\sqrt{MM^{\top}}$, where the
superscript $\top$ denotes the transpose of vectors or matrices. We say
$M>($resp. $\geq)\ 0$ if $M\in\mathbb{S}^{n}$ is positive (resp. nonnegative) definite.

Let $W(\cdot)=(W_{1}(\cdot),\ldots,W_{d}(\cdot))^{\top}$ be a standard
$d$-dimensional Brownian motion defined on a complete probability space
$(\Omega,\mathcal{F},P)$. The information structure is given by a filtration
$F=\{\mathcal{F}_{t}\}_{0\leq t\leq T},$ which is generated by $W(\cdot)$ and
augmented by all the $P$-null sets. Let $L_{\mathcal{F}_{T}}^{2}(\Omega;H)$
(resp. $L_{\mathcal{F}_{T}}^{\infty}(\Omega;H)$) denote the space of all
$\mathcal{F}_{T}$- measurable random variable $\eta$ with values in $H$ such
that%
\[
\mathbb{E}[|\eta|^{2}]<\infty\text{ (resp. }\eta\text{ is uniformly
bounded)};
\]
let $L_{\mathcal{F}}^{2}(0,T;H)$ (resp. $L_{\mathcal{F}}^{\infty}(0,T;H)$)
denote the space of all $\mathcal{F}_{t}$-progressively measurable processes
$x(\cdot)$ with values in $H$ such that%
\[
\mathbb{E}[\int_{0}^{T}|x(t)|^{2}dt]<\infty\text{ (resp. }\underset{(t,\omega
)\in\lbrack0,T]\times\Omega}{ess\sup}\mid x(t,\omega)\mid<\infty\text{)}.
\]

To simplify the presentation, we assume the dimension of the Brownian motion
$d=1$.

For a given $x\in\mathbb{R}^{n}$, consider the following linear stochastic
differential equation:
\begin{equation}
\left\{
\begin{array}
[c]{rcl}%
dX(t) & = & [A(t)X(t)+B(t)u(t)+b(t)]dt+[C(t)X(t)+D(t)u(t)+\sigma
(t)]dW_{t}\ ,\ t\in\lbrack0,T],\\
X(0) & = & x,
\end{array}
\right.  \label{state-equ}%
\end{equation}
where $A,B,C,D,b,\sigma$ are deterministic matrix-valued functions of suitable
sizes. In the above equation, $u(\cdot)$ is a control process and $X(\cdot)$
is the corresponding state process. In addition, the quadratic cost functional
is given by
\begin{equation}%
\begin{array}
[c]{rcl}%
J(u(\cdot)) & = & \mathbb{E}\{\frac{1}{2}\int_{0}^{T}%
[<Q(t)X(t),X(t)>+2<S(t)X(t),u(t)>+<R(t)u(t),u(t)>]dt\\
&  & +\frac{1}{2}<GX(T),X(T)>\},
\end{array}
\label{object}%
\end{equation}
where $G\in\mathbb{S}^{n}$, $Q,$ $S$ and $R$ are $\mathbb{S}^{n}-$,
$\mathbb{R}^{k\times n}-$ and $\mathbb{S}^{k}-$valued functions, respectively.

We introduce the following assumption:

\begin{assumption}
\label{assum2.1} Suppose that
\[
\left\{
\begin{array}
[c]{lc}%
A,C\in L^{\infty}(0,T;\mathbb{R}^{n\times n}),\ B,D,S^{\top}\in L^{\infty
}(0,T;\mathbb{R}^{n\times k}); & \\
Q\in L^{\infty}(0,T;\mathbb{S}^{n}),\ R\in L^{\infty}(0,T;\mathbb{S}%
^{k}),\ G\in\mathbb{S}^{n}; & \\
b,\sigma\in L^{2}(0,T;\mathbb{R}^{n}). &
\end{array}
\right.
\]

\end{assumption}

Let $B=\{0,1\}^{k}$, $C\subseteq\mathbb{R}^{k}$ be a closed and convex set,
$U=C\cap B\neq\varnothing$, $e=(1,1...1)^{\top}\in\mathbb{R}^{k}$ and
$E=\{u\in\mathbb{R}^{k}:0\leq u\leq e\}$. We set
\[
\mathcal{U}_{ad}=\{u(\cdot)\in L_{\mathcal{F}}^{2}(0,T;\mathbb{R}^{k})|u(t)\in
U\ \ a.e.\ a.s.\}.
\]
An element of $\mathcal{U}_{ad}$ is called an admissible control. Under
Assumption \ref{assum2.1}, there exists a unique solution $X(\cdot)\in
L_{\mathcal{F}}^{2}(0,T;\mathbb{R}^{n})$ to ($\ref{state-equ})$, which is
called the corresponding admissible state process. $(X(\cdot),u(\cdot))$ is
called an admissible pair.

Our stochastic linear quadratic control problem is to find an admissible
control $\bar{u}(t)$ such that
\begin{equation}
J(\bar{u}(\cdot))=\min\limits_{u(\cdot)\in\mathcal{U}_{ad}}J(u(\cdot
)).\label{opti object}%
\end{equation}

\section{Results of Linear Quadratic Problem by Functional Analysis}

Since the state equation in a stochastic linear quadratic problem is linear,
by variation of constant formula, the state process can be explicitly
expressed in terms of the initial state and the control. Substituting this
relation into the cost functional, we obtain a functional quadratic in the
state and control terms. To describe the method in detail, we first introduce
the following matrix-valued process:
\[
\left\{
\begin{array}
[c]{rcl}%
d\Phi(t) & = & A(t)\Phi(t)dt+C(t)\Phi(t)dW_{t},\ t\geq0,\\
\Phi(0) & = & I,
\end{array}
\right.
\]
from the method of stochastic differential equation, we know that $\Phi
^{-1}(t)$ exists for all $t\geq0$, and satisfy
\[
\left\{
\begin{array}
[c]{rcl}%
d\Phi^{-1}(t) & = & -\Phi^{-1}(t)[A(t)-C^{2}(t)]dt-\Phi^{-1}(t)C(t)dW_{t}%
,\ \ t\geq0,\\
\Phi^{-1}(0) & = & I.
\end{array}
\right.
\]
The solution of (\ref{state-equ}) can be written as:
\[%
\begin{array}
[c]{rcl}%
X(t) & = & \Phi(t)x+\Phi(t)\int_{0}^{t}\Phi^{-1}%
(s)[(B(s)-C(s)D(s))u(s)+b(s)-C(s)\sigma(s)]ds\\
&  & +\Phi(t)\int_{0}^{t}\Phi^{-1}(t)[D(s)u(s)+\sigma(s)]dW_{s},\ \ t\in
\lbrack0,T].
\end{array}
\]
By BDG-inequality, we can get the estimation of $X(t)$ as
\[
\mathbb{E}[\sup\limits_{s\in\lbrack0,t]}|X(s)|^{2}]\leq K\mathbb{E}%
\{|x|^{2}+\int_{0}^{t}[|u(s)|^{2}+|b(s)|^{2}+|\sigma(s)|^{2}]ds\},\ \ \forall
t\in\lbrack0,T].
\]
Next, $\forall x\in\mathbb{R}^{n}$, $u(\cdot)\in\mathcal{U}_{ad}$, we define
the following operators:
\[
\left\{
\begin{array}
[c]{l}%
(Lu(\cdot))(\cdot):=\Phi(\cdot)\{\int_{0}^{\cdot}\Phi^{-1}%
(s)[B(s)-C(s)D(s)]u(s)ds+\int_{0}^{\cdot}\Phi^{-1}(s)D(s)u(s)dW_{s}\},\\
\hat{L}u(\cdot):=(Lu(\cdot))(T),\ \ \forall u(\cdot)\in\mathcal{U}_{ad},\\
(\Gamma x)(\cdot):=\Phi(\cdot)x,\ \ \hat{\Gamma}x:=(\Gamma x)(T),\ \ \forall
x\in\mathbb{R}^{n},\\
f(\cdot):=\Phi(\cdot)\{\int_{0}^{\cdot}\Phi^{-1}(s)[b(s)-C(s)\sigma
(s)]ds+\int_{0}^{\cdot}\Phi^{-1}(s)\sigma(s)dW_{s}\},\\
\hat{f}:=f(T).
\end{array}
\right.
\]
Then, the state equation (\ref{state-equ}) and its terminal value can be
written as:
\[
\left\{
\begin{array}
[c]{l}%
X(\cdot)=(\Gamma x)(\cdot)+(Lu(\cdot))(\cdot)+f(\cdot),\\
X(T)\ =\hat{\Gamma}x+\hat{L}u(\cdot)+\hat{f}.
\end{array}
\right.
\]
Our next goal is to find a representation of the cost functional
(\ref{object}) in terms of control. To this end, we note that the following
operators are bounded linear operators
\[
\left\{
\begin{array}
[c]{l}%
L:\mathcal{U}_{ad}\rightarrow L_{\mathcal{F}}^{2}(0,T;\mathbb{R}^{n}%
),\ \hat{L}:\mathcal{U}_{ad}\rightarrow L_{\mathcal{F}_{T}}^{2}(\Omega
;\mathbb{R}^{n}),\\
\Gamma:\mathbb{R}^{n}\rightarrow L_{\mathcal{F}}^{2}(0,T;\mathbb{R}%
^{n}),\ \hat{\Gamma}:\mathbb{R}^{n}\rightarrow L_{\mathcal{F}_{T}}^{2}%
(\Omega;\mathbb{R}^{n}).
\end{array}
\right.
\]
We need to find the adjoint operators of the above bounded linear operators,
\[
\left\{
\begin{array}
[c]{l}%
L^{\ast}:L_{\mathcal{F}}^{2}(0,T;\mathbb{R}^{n})\rightarrow\mathcal{U}%
_{ad},\ \hat{L}^{\ast}:L_{\mathcal{F}_{T}}^{2}(\Omega;\mathbb{R}%
^{n})\rightarrow\mathcal{U}_{ad},\\
\Gamma^{\ast}:L_{\mathcal{F}}^{2}(0,T;\mathbb{R}^{n})\rightarrow\mathbb{R}%
^{n},\ \hat{\Gamma}^{\ast}:L_{\mathcal{F}_{T}}^{2}(\Omega;\mathbb{R}%
^{n})\rightarrow\mathbb{R}^{n},
\end{array}
\right.
\]
such that:
\begin{equation}
\left\{
\begin{array}
[c]{l}%
\mathbb{E}[\int_{0}^{T}<(Lu(\cdot))(t),\xi(t)>dt]=\mathbb{E}[\int_{0}%
^{T}<u(\cdot),(L^{\ast}\xi(\cdot))(t)>dt],\\
\mathbb{E}[\int_{0}^{T}<(\Gamma x)(t),\xi(t)>dt]=\mathbb{E}[<x,(\Gamma^{\ast
}\xi(\cdot))>],\\
\ \ \ \forall u(\cdot)\in\mathcal{U}_{ad},\ \xi(\cdot)\in L_{\mathcal{F}}%
^{2}(0,T;\mathbb{R}^{n}),
\end{array}
\right.  \label{L*}%
\end{equation}
and
\begin{equation}
\left\{
\begin{array}
[c]{l}%
\mathbb{E}[<\hat{L}u(\cdot),\eta>]=\mathbb{E}[\int_{0}^{T}<u(t),(\hat{L}%
^{\ast}\eta)(t)>dt],\\
\mathbb{E}[<\hat{\Gamma}x,\eta>]=\mathbb{E}[<x,(\hat{\Gamma}^{\ast}\eta>],\\
\ \ \ \forall u(\cdot)\in\mathcal{U}_{ad},\ \eta\in L_{\mathcal{F}}^{2}%
(\Omega;\mathbb{R}^{n}),\ \ x\in\mathbb{R}^{n}.
\end{array}
\right.  \label{L-hat}%
\end{equation}
Actually, we can define the adjoint operators through the following backward
stochastic differential equation:
\begin{equation}
\left\{
\begin{array}
[c]{ll}%
dp(t)=-[A(t)^{\top}p(t)+C(t)^{\top}q(t)+\xi(t)]dt+q(t)dW_{t}, & \ \ t\in
\lbrack0,T],\\
p(T)=\eta\in L_{\mathcal{F}}^{2}(\Omega;\mathbb{R}^{n}). &
\end{array}
\right.  \label{adjoint1}%
\end{equation}
We have the following results. \newline

\begin{proposition}
(\cite{yongzhou}) \label{propo3.1} $(i)$ For any $\xi(\cdot)\in L_{\mathcal{F}%
}^{2}(0,T;\mathbb{R}^{n})$, let $(p_{0}(\cdot),q_{0}(\cdot))\in L_{\mathcal{F}%
}^{2}(0,T;\mathbb{R}^{n})\times L_{\mathcal{F}}^{2}(0,T;\mathbb{R}^{n})$ be
the adapted solution of $(\ref{adjoint1})$ with $\eta=0$. Define
\[
\left\{
\begin{array}
[c]{ll}%
(L^{\ast}\xi)(t):=B(t)^{\top}p_{0}(t)+D(t)^{\top}q_{0}(t), & \ \ t\in\lbrack0,T],\\
\Gamma^{\ast}\xi=p_{0}(0). &
\end{array}
\right.
\]
Then $L^{\ast}:L_{\mathcal{F}}^{2}(0,T;\mathbb{R}^{n})\rightarrow
\mathcal{U}_{ad}$ and $\Gamma^{\ast}:L_{\mathcal{F}}^{2}(0,T;\mathbb{R}%
^{n})\rightarrow\mathbb{R}^{n}$ are bounded operators satisfying $(\ref{L*})$.

$(ii)$ For any $\eta\in L_{\mathcal{F}}^{2}(\Omega;\mathbb{R}^{n})$, let
$(p_{1}(\cdot),q_{1}(\cdot))\in L_{\mathcal{F}}^{2}(0,T;\mathbb{R}^{n})\times
L_{\mathcal{F}}^{2}(0,T;\mathbb{R}^{n})$ be the adapted solution of
$(\ref{adjoint1})$ with $\xi(\cdot)=0$. Define
\[
\left\{
\begin{array}
[c]{ll}%
(\hat{L}^{\ast}\eta)(t):=B(t)^{\top}p_{1}(t)+D(t)^{\top}q_{1}(t), & \ \ t\in
\lbrack0,T],\\
\hat{\Gamma}^{\ast}\eta=p_{1}(0). &
\end{array}
\right.
\]
Then $\hat{L}^{\ast}:L_{\mathcal{F}_{T}}^{2}(\Omega;\mathbb{R}^{n}%
)\rightarrow\mathcal{U}_{ad}$ and $\hat{\Gamma}^{\ast}:L_{\mathcal{F}_{T}}%
^{2}(\Omega;\mathbb{R}^{n})\rightarrow\mathbb{R}^{n}$ are bounded operators
satisfying $(\ref{L-hat})$.
\end{proposition}

Once we have the above results, we can obtain another form of the cost
functional:
\begin{equation}
J(u(\cdot))=\frac{1}{2}\{<Nu(\cdot),u(\cdot)>+2<H(x),u(\cdot
)>+M(x)\},\label{functional-object}%
\end{equation}
where
\begin{equation}
\left\{
\begin{array}
[c]{l}%
N=R+L^{\ast}QL+SL+L^{\ast}S^{\top}+\hat{L}^{\ast}G\hat{L},\\
H(x)=(L^{\ast}Q+S)[(\Gamma x)(\cdot)+f(\cdot)]+\hat{L}^{\ast}G(\hat{\Gamma
}x+\hat{f}),\\
M(x)=<Q[(\Gamma x)(\cdot)+f(\cdot)],(\Gamma x)(\cdot)+f(\cdot)>\\
\ \ \ \ \ \ \ \ \ \ \ \ +<G(\hat{\Gamma}x+\hat{f}),\hat{\Gamma}x+\hat{f}>.
\end{array}
\right.  \label{functional-coeffi}%
\end{equation}
It is clear that the operator $N:\mathcal{U}_{ad}\rightarrow\mathcal{U}_{ad}$
is a bounded linear operator.

\section{Stochastic Maximum Principle}

Generally speaking, since the control domain of the stochastic linear
quadratic problem is nonconvex, one should take the second-order adjoint
process into consideration. To conveniently state this classical maximum
principle, we give the following results, which come from \cite{yongzhou}
Chapter 3. We consider the following controlled stochastic differential
equation:
\begin{equation}
\left\{
\begin{array}
[c]{ll}%
dX(t)=b(t,X(t),u(t))dt+\sigma(t,X(t),u(t))dW_{t}, & \ \ t\in\lbrack0,T],\\
X(0)=x, &
\end{array}
\right.  \label{nonconv-syst}%
\end{equation}
with the cost functional
\begin{equation}
J(u(\cdot))=\mathbb{E}[\int_{0}^{T}f(t,X(t),u(t))dt+h(X(T))].
\label{nonconv-object}%
\end{equation}
The controller wants to find the infimum of $J$ over admissible control set.
Moreover, we introduce the following assumption as in \cite{yongzhou}.

\begin{assumption}
\label{assum4.1} $(i)$ $(U,d)$ is a separable metric space and $T>0$%
.\newline$(ii)$ The maps $b,\sigma,f$, and $h$ are measurable, and there exist
a constant $L>0$ and a modulus of continuity $\bar{\omega}:[0,\infty
)\rightarrow\lbrack0,\infty)$ such that for $\phi(t,x,u)=b(t,x,u),\sigma
(t,x,u),f(t,x,u)$, $h(x)$, we have
\[
\left\{
\begin{array}
[c]{ll}%
|\phi(t,x,u)-\phi(t,\hat{x},\hat{u})| & \leq L|x-\hat{x}|+\bar{\omega
}(d(u,\hat{u})),\\
& \forall t\in\lbrack0,T],x,\hat{x}\in\mathbb{R}^{n},u,\hat{u}\in U,\\
|\phi(t,0,u)|\leq L, & \forall(t,u)\in\lbrack0,T]\times U.
\end{array}
\right.
\]
$(iii)$ The maps $b,\sigma,f$, and $h$ are $C^{2}$ in $x$. Moreover, there
exist a constant $L>0$ and a modulus of continuity $\bar{\omega}%
:[0,\infty)\rightarrow\lbrack0,\infty)$ such that for $\phi
(t,x,u)=b(t,x,u),\sigma(t,x,u)$, $f(t,x,u),h(x)$, we have
\[
\left\{
\begin{array}
[c]{ll}%
|\phi_{x}(t,x,u)-\phi_{x}(t,\hat{x},\hat{u})| & \leq L|x-\hat{x}|+\bar{\omega
}(d(u,\hat{u})),\\
|\phi_{xx}(t,x,u)-\phi_{xx}(t,\hat{x},\hat{u})| & \leq\bar{\omega}(|x-\hat
{x}|+d(u,\hat{u})),\\
& \forall t\in\lbrack0,T],x,\hat{x}\in\mathbb{R}^{n},u,\hat{u}\in U.
\end{array}
\right.
\]

\end{assumption}

Let $\bar{X}(\cdot),\bar{u}(\cdot)$ be an optimal pair of the stochastic
system (\ref{nonconv-syst})-(\ref{nonconv-object}). We introduce the following
first-and second-order adjoint equation
\begin{equation}
\left\{
\begin{array}
[c]{l}%
dp(t)=-[b_{x}(t,\bar{X}(t),\bar{u}(t))^{\top}p(t)+\sigma_{x}(t,\bar{X}(t),\bar
{u}(t))^{\top}q(t)-f_{x}(t,\bar{X}(t),\bar{u}(t))]dt+q(t)dW_{t},\ t\in
\lbrack0,\ T],\\
p(T)=-h_{x}(\bar{X}(T)),
\end{array}
\right.  \label{first-adjoint}%
\end{equation}

\begin{equation}
\left\{
\begin{array}
[c]{lcl}%
dP(t) & = & -[b_{x}(t,\bar{X}(t),\bar{u}(t))^{\top}P(t)+P(t)b_{x}(t,\bar
{X}(t),\bar{u}(t))\\
&  & +\sigma_{x}(t,\bar{X}(t),\bar{u}(t))^{\top}P(t)\sigma_{x}(t,\bar{X}%
(t),\bar{u}(t))\\
&  & +\sigma_{x}(t,\bar{X}(t),\bar{u}(t))^{\top}Q(t)+Q(t)\sigma_{x}(t,\bar
{X}(t),\bar{u}(t))\\
&  & +H_{xx}(t,\bar{X}(t),\bar{u}(t),p(t),q(t))]dt+Q(t)dW_{t},\\
P(T) & = & -h_{xx}(\bar{X}(T)),
\end{array}
\right.  \label{second-adjoint}%
\end{equation}
where the Hamiltonian $H$ is defined by
\[%
\begin{array}
[c]{lcl}%
H(t,x,u,p,q) & = & <p,b(t,x,u)>+q^{\top}\sigma(t,x,u)-f(t,x,u),\\
&  & (t,x,u,p,q)\in\lbrack0,T]\times\mathbb{R}^{n}\times U\times\mathbb{R}%
^{n}\times\mathbb{R}^{n}.
\end{array}
\]

\begin{proposition}
\label{propo4.2} (General Stochastic Maximum principle \cite{peng1990}
\cite{yongzhou}) Let Assumption \ref{assum4.1} hold, $(\bar{X}(\cdot),\bar
{u}(\cdot))$ be an optimal pair of the stochastic system (\ref{nonconv-syst}%
)-(\ref{nonconv-object}). Then there exist pairs of process
\[
\left\{
\begin{array}
[c]{l}%
(p(\cdot),q(\cdot))\in L_{\mathcal{F}}^{2}(0,T;\mathbb{R}^{n})\times
L_{\mathcal{F}}^{2}(0,T;\mathbb{R}^{n}),\\
(P(\cdot),Q(\cdot))\in L_{\mathcal{F}}^{2}(0,T;\mathbb{S}^{n})\times
L_{\mathcal{F}}^{2}(0,T;\mathbb{S}^{n}),
\end{array}
\right.
\]
satisfying the first-and second-order adjoint equations $(\ref{first-adjoint}%
)$ and $(\ref{second-adjoint})$, respectively, such that
\[%
\begin{array}
[c]{l}%
H(t,\bar{X}(t),\bar{u}(t),p(t),q(t))-H(t,\bar{X}(t),v,p(t),q(t))\\
\ \ -\frac{1}{2}(\sigma(t,\bar{X}(t),\bar{u}(t))-\sigma(t,\bar{X}%
(t),v))^{\top}P(t)(\sigma(t,\bar{X}(t),\bar{u}(t))-\sigma(t,\bar{X}(t),v))\geq0,\\
\ \ \ \ \ \ \ \ \ \ \ \ \ \ \ \ \ \ \ \ \ \ \ \ \ \ \ \ \ \ \ \ \ \ \ \ \ \ \ \ \ \ \ \ \ \ \ \ \ \ \forall
v\in U,a.e.\ t\in\lbrack0,T],\mathbb{P}-a.s..
\end{array}
\]

\end{proposition}

It should be noted that the second-order adjoint equations have to be
introduced in the above stochastic maximum principles. Different from
\cite{peng1990} \cite{yongzhou}, the primal problem will be turned into a
concave problem with convex control domain in this section. \ And a novel
stochastic maximum principle without the second-order adjoint equation will be also
proposed. To proceed further, the following useful Proposition is needed. The
interested readers can find this in \cite{peng1993} and \cite{yongzhou}.

\begin{proposition}
\label{propo4.3} (Stochastic Maximum principle \cite{peng1993}) Let $\bar
{u}(\cdot)$ be an optimal control and let $\bar{X}(\cdot)$ be the
corresponding trajectory. The control domain is convex and all the
coefficients are $C^{1}$ in $u$, we have
\[%
\begin{array}
[c]{l}%
<H_{u}(t,\bar{X}(t),\bar{u}(t),p(t),q(t)),v-\bar{u}(t)>\ \leq\ 0,\\
\ \ \ \ \ \ \ \ \ \ \ \ \ \ \ \ \forall v\in U,a.e.\ t\in\lbrack
0,T],\mathbb{P}-a.s..
\end{array}
\]

\end{proposition}

The next theorem is one of the main results in this literature. We first transform
the original linear quadratic problem into a quadratic optimization problem in
a Hilbert space by functional analysis approach as shown in Section 3. Then
through the introduced parameter, we can turn the original problem into a
concave control problem with convex control domain. Then we can apply the
stochastic maximum principle in Proposition \ref{propo4.3} to the transformed
concave problem, and obtain the new stochastic maximum principle.

\begin{theorem}
\label{theorem} Let Assumption \ref{assum2.1} hold, and let $(\bar{X}%
(\cdot),\bar{u}(\cdot))$ be an optimal pair of the stochastic linear quadratic
control problem $(\ref{state-equ})-(\ref{opti object})$. Then there exists an
adapted solution $(p(\cdot),q(\cdot))\in L_{\mathcal{F}}^{2}(0,T;\mathbb{R}%
^{n})\times L_{\mathcal{F}}^{2}(0,T;\mathbb{R}^{n})$ to the following backward
stochastic differential equation
\[
\left\{
\begin{array}
[c]{l}%
dp(t)=-[A(t)^{\top}p(t)+C(t)^{\top}q(t)-Q(t)\bar{X}(t)-S(t)^{\top}\bar{u}%
(t)]dt+q(t)dW_{t},\\
p(T)=-G\bar{X}(T),
\end{array}
\right.
\]
such that
\[%
\begin{array}
[c]{r}%
<H_{u}^{\mu}(t,\bar{X}(t),\bar{u}(t),p(t),q(t)),v-\bar{u}(t)>\ \leq0,\\
\forall v\in\bar{U},a.e.\ t\in\lbrack0,T],\mathbb{P}-a.s.,
\end{array}
\]
where the Hamiltonian function $H^{\mu}$ is defined by
\begin{equation}%
\begin{array}
[c]{lll}%
H^{\mu}(t,x,u,p,q) & = & <p,A(t)x+B(t)u+b(t)>+<q,C(t)x+D(t)u+\sigma(t)>\\
&  & -\frac{1}{2}\{[<Q(t)x,x>+2<S(t)x-\mu\mathbf{I},u>+<(R(t)+\mu
\mathbf{I})u(t),u(t)>]\},
\end{array}
\label{new-hamil}%
\end{equation}
where the parameter $\mu$ is defined as: $-\mu$ is the largest eigenvalue of
$N$, which is described in $(\ref{functional-coeffi})$; the operator
$\mathbf{I}$ is identity operator.
\end{theorem}

\textbf{Proof.} The proof is technical. We divide it into two steps to make
the idea clear.

\emph{Step 1. \ \ Equivalent Formulation} \newline As shown in section 3, we
obtain a functional quadratic in the triple of the initial state, the control
and the nonhomogeneous term. That is
\[
J(u(\cdot))=\frac{1}{2}\{<Nu(\cdot),u(\cdot)>+2<H(x),u(\cdot)>+M(x)\},
\]
where
\[
\left\{
\begin{array}
[c]{l}%
N=R+L^{\ast}QL+SL+L^{\ast}S^{\top}+\hat{L}^{\ast}G\hat{L},\\
H(x)=(L^{\ast}Q+S)[(\Gamma x)(\cdot)+f(\cdot)]+\hat{L}^{\ast}G(\hat{\Gamma
}x+\hat{f}),\\
M(x)=<Q[(\Gamma x)(\cdot)+f(\cdot)],(\Gamma x)(\cdot)+f(\cdot)>\\
\ \ \ \ \ \ \ \ \ \ \ \ +<G(\hat{\Gamma}x+\hat{f}),\hat{\Gamma}x+\hat{f}>.
\end{array}
\right.
\]
Given any real number $\mu\in\mathbb{R}$, let $\bar{N}=N+\mu\mathbf{I}$ and
$\bar{H}(x)=H(x)-\frac{1}{2}\mu\mathbf{I}$, then we can define the following
stochastic linear quadratic control problem
\[%
\begin{array}
[c]{rcl}%
J^{\mu}(u(\cdot)) & = & \frac{1}{2}\{<\bar{N}u(\cdot),u(\cdot)>+2<\bar
{H}(x),u(\cdot)>+M(x)\}\\
& = & \frac{1}{2}\{<(N+\mu\mathbf{I})u(\cdot),u(\cdot)>+2<H(x)-\frac{1}{2}%
\mu\mathbf{I},u(\cdot)>+M(x)\}\\
& = & \mathbb{E}\{\frac{1}{2}\int_{0}^{T}[<Q(t)X(t),X(t)>+2<S(t)X(t)-\frac
{1}{2}\mu\mathbf{I},u(t)>\\
&  & +<(R(t)+\mu\mathbf{I})u(t),u(t)>]dt+\frac{1}{2}<GX(T),X(T)>\}.
\end{array}
\]
When we select $-\mu$ as the largest eigenvalue of $N$, then the operator
$\bar{N}$ is negative semi-definite, then $J^{\mu}$ is concave with respect to
$u(\cdot)$. We claim that the problem
\begin{equation}
\min\limits_{u(\cdot)\in\mathcal{U}_{ad}}J(u(\cdot)) \label{equiva1}%
\end{equation}
is equivalent (in the sense of coinciding optimal control and optimal value)
to
\begin{equation}
\min\limits_{u(\cdot)\in\bar{\mathcal{U}}_{ad}}J^{\mu}(u(\cdot)),
\label{equiva2}%
\end{equation}
where $\bar{\mathcal{U}}_{ad}$ is defined by
\[
\bar{\mathcal{U}}_{ad}=\{u(\cdot)\in L_{\mathcal{F}}^{2}(0,T;\mathbb{R}%
^{k})|u(t)\in\bar{U}=C\cap E\ \ a.e.\ a.s.\}.
\]
Note that we have defined $E=\{u\in\mathbb{R}^{k}:0\leq u\leq e\}$ in section 2.

Indeed, since $J^{\mu}(u(\cdot))$ is concave with respect to $u$, the optimal
control must attained on a vertex of $E$. Then $\min\limits_{u(\cdot
)\in\mathcal{U}_{ad}}J(u(\cdot))=\min\limits_{u(\cdot)\in\bar{\mathcal{U}%
}_{ad}}J^{\mu}(u(\cdot))$. If $(\bar{X}(\cdot),\bar{u}(\cdot))$ is the optimal
pair of problem $(\ref{equiva1})$, then $\forall u(\cdot)\in\bar{\mathcal{U}%
}_{ad}$,
\[%
\begin{array}
[c]{rcl}%
J^{\mu}(u(\cdot)) & = & \mathbb{E}\{\frac{1}{2}\int_{0}^{T}%
[<Q(t)X(t),X(t)>+2<S(t)X(t)-\frac{1}{2}\mu\mathbf{I},u(t)>\\
&  & +<(R(t)+\mu\mathbf{I})u(t),u(t)>]dt+\frac{1}{2}<GX(T),X(T)>\}\\
& \geq & \min\limits_{u(\cdot)\in\bar{\mathcal{U}}_{ad}}J^{\mu}(u(\cdot))\\
& = & \min\limits_{u(\cdot)\in\mathcal{U}_{ad}}J(u(\cdot))\\
& = & J(\bar{u}(\cdot))\\
& = & J^{\mu}(\bar{u}(\cdot)).
\end{array}
\]
That is $\bar{u}$ is the optimal control of $J^{\mu}(u(\cdot))$. \newline On
the other hand, if $(\bar{X}(\cdot),\bar{u}(\cdot))$ is the optimal pair of
problem $(\ref{equiva2})$, then $\forall u(\cdot)\in\mathcal{U}_{ad}$,
\[%
\begin{array}
[c]{rcl}%
J(u(\cdot)) & = & \mathbb{E}\{\frac{1}{2}\int_{0}^{T}%
[<Q(t)X(t),X(t)>+2<S(t)X(t),u(t)>+<R(t)u(t),u(t)>]dt\\
&  & +\frac{1}{2}<GX(T),X(T)>\}\\
& \geq & \min\limits_{u(\cdot)\in\mathcal{U}_{ad}}J(u(\cdot))\\
& = & \min\limits_{u(\cdot)\in\bar{\mathcal{U}}_{ad}}J^{\mu}(u(\cdot))\\
& = & J^{\mu}(\bar{u}(\cdot))\\
& = & J(\bar{u}(\cdot)).
\end{array}
\]
That is $\bar{u}$ is the optimal control of $J(u(\cdot))$. Therefore, the
problem $(\ref{equiva1})$ is equivalent to the minimization of a concave
quadratic function over the convex set, so we can apply the classical
stochastic maximum principle as shown in Proposition \ref{propo4.3}, which
deals with convex control domain, to the equivalent problem.

\emph{Step 2. \ \ Apply Stochastic Maximum Principle} \newline In order to
derive the maximum principle for the classical stochastic optimal control
problem, one needs to obtain the variational equation of state equation. Based
on the variational inequality and the introduced adjoint equation, one can get
the stochastic maximum principle (see \cite{ben1982},\cite{ben1983},
\cite{peng1993} and \cite{yongzhou}).

In our case, we should introduce the adjoint equation as the following
backward stochastic differential equation:
\[
\left\{
\begin{array}
[c]{l}%
dp(t)=-[A(t)^{\top}p(t)+C(t)^{\top}q(t)-Q(t)\bar{X}(t)-S(t)^{\top}\bar{u}%
(t)]dt+q(t)dW_{t},\\
p(T)=-G\bar{X}(T).
\end{array}
\right.
\]
Moreover, the Hamiltonian function $H^{\mu}$ is defined by
\[%
\begin{array}
[c]{lll}%
H^{\mu}(t,x,u,p,q) & = & <p,A(t)x+B(t)u+b(t)>+<q,C(t)x+D(t)u+\sigma(t)>\\
&  & -\frac{1}{2}\{[<Q(t)x,x>+2<S(t)x-\frac{\mu}{2}\mathbf{I},u>+<(R(t)+\mu
\mathbf{I})u(t),u(t)>]\},
\end{array}
\]
then, we can obtain the following stochastic maximum principle
\[%
\begin{array}
[c]{r}%
<H_{u}^{\mu}(t,\bar{X}(t),\bar{u}(t),p(t),q(t)),v-\bar{u}(t)>\ \leq0,\\
\forall v\in\bar{U},a.e.\ t\in\lbrack0,T],\mathbb{P}-a.s..
\end{array}
\]
This completes the proof. \hfill$\Box$ \smallskip\vskip4mm

\begin{remark}
\label{remark1} From the above theorem, we can get some further subtle
results. We can take into account the optimal control respectively. Then from
\[%
\begin{array}
[c]{rll}%
H_{u}^{\mu}(t,\bar{X}(t),\bar{u}(t),p(t),q(t)) & = & B(t)^{\top}p(t)+D(t)^{\top}%
q(t)\\
&  & -\frac{1}{2}[2S(t)\bar{X}(t)+(R(t)+\mu\mathbf{I})^{\top}\bar{u}%
(t)+(R(t)+\mu\mathbf{I})\bar{u}(t)^{\top}-\mu\mathbf{I}],
\end{array}
\]
we get

$\ \ \ \ H_{u}^{\mu}(t,\bar{X}(t),\bar{u}(t),p(t),q(t))=\left\{
\begin{array}
[c]{ll}%
B(t)^{\top}p(t)+D(t)^{\top}q(t)-S(t)\bar{X}(t)+\mu I\leq0,\  & if\ \ \bar{u}(t)=0;\\
B(t)^{\top}p(t)+D(t)^{\top}q(t)-R(t)-S(t)\bar{X}(t)-\frac{\mu}{2}\mathbf{I}\geq0, &
\ if\ \ \bar{u}(t)=e.
\end{array}
\right.  $

We can compare the above results with the primal maximum principle which deals
with nonconvex control domain problem as shown in Proposition \ref{propo4.2}.
Because the control domain is nonconvex, one should introduce first-and
second-order adjoint equations as
\[
\left\{
\begin{array}
[c]{l}%
dp(t)=-[A(t)^{\top}p(t)+C(t)^{\top}q(t)-Q(t)\bar{X}(t)-S(t)^{\top}\bar{u}%
(t)]dt+q(t)dW_{t},\\
p(T)=-G\bar{X}(T),
\end{array}
\right.
\]
and
\[
\left\{
\begin{array}
[c]{lll}%
dP(t) & = & -[A(t)^{\top}P(t)+P(t)A(t)+C(t)^{\top}P(t)C(t)+\Lambda(t)C(t)+C(t)^{\top}%
\Lambda(t)-Q(t)]dt\\
&  & +\Lambda(t)dW_{t},\\
P(T) & = & -G.
\end{array}
\right.
\]
The maximum principle is
\[%
\begin{array}
[c]{rl}%
H^{0}(t,\bar{X}(t),\bar{u}(t),p(t),q(t))-H^{0}(t,\bar{X}(t),v,p(t),q(t)) & \\
-\frac{1}{2}(\bar{u}(t)-v)^{\top}D(t)^{\top}P(t)D(t)(\bar{u}(t)-v)\geq 0, & \\
\forall v\in U,a.e.\ t\in\lbrack0,T],\mathbb{P}-a.s., &
\end{array}
\]
where $H^{0}$ is defined by $(\ref{new-hamil})$ when $\mu=0$. Moreover, for
any $v\in U$, we have the following results
\[
\left\{
\begin{array}
[c]{ll}%
\frac{1}{2}v^{\top}[R(t)-D(t)^{\top}P(t)D(t)]v-v^{\top}[B(t)^{\top}p(t)+D(t)^{\top}%
q(t)-S(t)\bar{X}(t)]\geq0,\ \ \  & if\ \ \bar{u}(t)=0;\\
(e-v)^{\top}[B(t)^{\top}p(t)+D(t)^{\top}q(t)-S(t)\bar{X}(t)]+v^{\top}R(t)v-e^{\top}R(t)e & \\
\ \ \ -\frac{1}{2}(e-v)^{\top}D(t)^{\top}P(t)D(t)(e-v)\geq0,\ \ \  & if\ \ \bar
{u}(t)=e.
\end{array}
\right.
\]
We can see that the nonconvex maximum principle is complex. However, the
approach provided in this paper is relatively concise. In the next section we
will give an example to show how the new maximum principle works.
\end{remark}

{}

\section{Example}

\setcounter{equation}{0} Let us consider the following controlled stochastic
differential equation ($n=1$):
\[
\left\{
\begin{array}
[c]{rcl}%
dX(t) & = & u(t)dW_{t}\ ,\ t\in\lbrack0,1],\\
X(0) & = & 0,
\end{array}
\right.
\]
with the control domain $U=\{0,1\}\cap\lbrack0,1]=\{0,1\}$. The cost
functional is defined by
\[
J(u(\cdot))=\mathbb{E}\{\int_{0}^{1}[X(t)^{2}-\frac{1}{2}u(t)^{2}%
]dt+X(1)^{2}\}.
\]
Substituting $X(t)=\int_{0}^{t}u(s)dW_{s}$ into the cost functional, we obtain
the following:
\[%
\begin{array}
[c]{rcl}%
J(u(\cdot)) & = & \mathbb{E}\{\int_{0}^{1}[X(t)^{2}-\frac{1}{2}u(t)^{2}%
]dt+X(1)^{2}\}\\
& = & \mathbb{E}\{\int_{0}^{1}[(\int_{0}^{t}u(s)dW_{s})^{2}-\frac{1}%
{2}u(t)^{2}]dt+(\int_{0}^{1}u(s)dW_{s})^{2}\}\\
& = & \int_{0}^{1}[\mathbb{E}(\int_{0}^{t}u(s)dW_{s})^{2}]dt-\frac{1}%
{2}\mathbb{E}[\int_{0}^{1}u(t)^{2}dt]+\mathbb{E}(\int_{0}^{1}u(s)dW_{s})^{2}\\
& = & \mathbb{E}[\int_{0}^{1}\int_{0}^{t}u(s)^{2}dsdt]+\frac{1}{2}%
\mathbb{E}[\int_{0}^{1}u(t)^{2}dt]\\
& = & \mathbb{E}[\int_{0}^{1}\int_{s}^{1}u(s)^{2}dtds]+\frac{1}{2}%
\mathbb{E}[\int_{0}^{1}u(t)^{2}dt]\\
& = & \mathbb{E}[\int_{0}^{1}(1-t)u(t)^{2}dt]+\frac{1}{2}\mathbb{E}[\int%
_{0}^{1}u(t)^{2}dt]\\
& = & \mathbb{E}[\int_{0}^{1}(\frac{3}{2}-t)u(t)^{2}dt].
\end{array}
\]
Hence, the optimal control is $\bar{u}(t)=0$, the corresponding optimal state
trajectory is $\bar{X}(t)=0$. The first- and second adjoint equations
associated with the optimal pair are
\begin{equation}
\left\{
\begin{array}
[c]{rcl}%
dp(t) & = & q(t)dW_{t},\\
p(1) & = & 0,
\end{array}
\right.  \label{ex-adjoint}%
\end{equation}%
\[
\left\{
\begin{array}
[c]{rcl}%
dP(t) & = & 2dt+Q(t)dW_{t},\\
P(1) & = & -2.
\end{array}
\right.
\]
Thus, by the uniqueness of the solution to backward stochastic differential
equation, we get $(p(t),q(t))=(0,0)$, $(P(t),Q(t))=(2t-4,0)$. The
corresponding Hamiltonian is
\[
H(t,\bar{X}(t),u,p(t),q(t))=\frac{1}{2}u^{2}+q(t)u=\frac{1}{2}u^{2}.
\]
This is a convex function in $u$, which does not attain a maximum at $\bar
{u}(t)=0$ for any $t\in\lbrack0,1]$. But the $\mathcal{H}$-function
\[
\mathcal{H}(t,\bar{X}(t),u):=H(t,\bar{X}(t),u,p(t),q(t))-\frac{1}{2}%
P(t)\bar{u}(t)+\frac{1}{2}P(t)(u-\bar{u}(t))^{2}=\frac{1}{2}(2t-3)u^{2},
\]
is concave with respect to $u$ for $\forall t\in\lbrack0,1]$, and $\bar
{u}(t)=0$ does maximize $\mathcal{H}$. Thus we see clearly how the
second-order adjoint process plays a role in turning the convex function
$u\longmapsto H(t,\bar{X}(t),u,p(t),q(t))$ into the concave one $u\longmapsto
\mathcal{H}(t,\bar{X}(t),u)$.

As to the new approach which provided in this paper, in order to cope with the
primal problem, we introduce a parameter $\mu$. In this example we let
$\mu=-3$, then
\begin{equation}%
\begin{array}
[c]{rcl}%
J^{\mu}(u(\cdot)) & = & J(u(\cdot))+\frac{1}{2}<(-3)u(\cdot),u(\cdot
)>-\frac{1}{2}<(-3)\mathbf{I},u(\cdot)>\\
& = & \mathbb{E}\{\int_{0}^{1}[X(t)^{2}-\frac{1}{2}u(t)^{2}-\frac{3}%
{2}u(t)^{2}+\frac{3}{2}u(t)]dt+X(1)^{2}\},
\end{array}
\end{equation}
with $\bar{U}=[0,1]$. As we have pointed out that the problem is equivalent to
primal problem in the sense of coinciding optimal solutions. The adjoint
equation is the same as (\ref{ex-adjoint}), and the corresponding Hamiltonian
is
\[
H^{\mu}(t,\bar{X}(t),u,p(t),q(t))=2u^{2}-\frac{3}{2}u.
\]
From classical stochastic maximum principle, we obtain that the following inequality holds
when the control domain is a convex set:
\[
H_{u}^{\mu}(t,\bar{X}(t),u,p(t),q(t))|_{u=\bar{u}(t)}\cdot(v-\bar{u}%
(t))\leq0,\ \forall v\in\lbrack0,1].
\]
That is
\[
(-\frac{3}{2})\cdot(v-\bar{u}(t))\leq0,\ \ \forall v\in\lbrack0,1].
\]
From the stochastic maximum principle, we know that $\bar{u}(t)=0$ is a
candidate optimal control.

We can say more about this result. In this case, the control domain of primal
problem is nonconvex. One needs to study both the first- and second-order
terms in the Taylor expansion of the spike variation and come up with a
stochastic maximum principle. The second-order adjoint process plays a role in
turning convex Hamiltonian $H$ into concave one $\mathcal{H}$, as we have shown
above. However, in this paper, without second-order adjoint process,
we can obtain a stochastic maximum principle involving a parameter $\mu$.
{}

\section{Conclusion}

In this article, we investigate a new stochastic maximum principle to deal
with the stochastic linear quadratic problem with nonconvex control domain. We
turn the original problem into a concave control problem with convex
control domain by employing the functional analysis approach and introducing
a parameter. Then we apply the classical stochastic maximum principle to the
transformed concave problem. Compared with the existing methods, the developed
stochastic maximum principle is easy to check.  Finally, we provide a work-out example
to illustrate the effectiveness of the proposed method in this paper.
Readers can have a try to solve the problem by utilizing the Hamilton-Jacobi-Bellman
approach or other methods. Although considering the forward state equation, we encourage the
readers to explore backward stochastic system or forward-backward stochastic
system, and to obtain results that are not yet known to us.  \newline%
\textbf{Acknowledgment.} The second author thanks Dr. Haodong Liu and Xiaomin
Shi for useful discussions related to this work. \newline

\end{document}